\documentclass[12pt]{article}
\title{Oration for Andrew Wiles}
\date{}
\begin{document}
 \maketitle

\subsection*{Fanfare}
We honour Andrew Wiles for his supreme contribution to number theory, a
contribution that has made him the world's most famous mathematician and
a beacon of inspiration for students of math; while solving Fermat's
Last Theorem, for 350 years the most celebrated open problem in
mathematics, Wiles's work has also dramatically opened up whole new
areas of research in number theory.

\subsection*{A love of mathematics}
The bulk of this eulogy is mathematical, for which I make no apology. I
want to stress here that, in addition to calculations in which each line
is correctly deduced from the preceding lines, mathematics is above all
passion and drama, obsession with solving the unsolvable. In a modest
way, many of us at Warwick share Andrew Wiles' overriding passion for
mathematics and its unsolved problems.

\subsection*{Three short obligatory pieces}
\paragraph{Biography}
 Oxford, Cambridge, Royal Society Professor at Oxford from 1988,
 Professor at Princeton since 1982 (lamentably for maths in Britain).
 Very many honours in the last 5 years, including the Wolf prize, Royal
 Society gold medal, the King Faisal prize, many, many others.
\paragraph{Human interest story}
 The joy and pain of Wiles's work on Fermat are beautifully documented
 in John Lynch's BBC Horizon documentary; I particularly like the bit
 where Andrew takes time off from unravelling the riddle that has
 baffled the world's best minds for 350 years to tell bed-time stories
 to little Clare, Kate and Olivia.
\paragraph{Predictable barbed comment on Research Assessment}
 It goes without saying that an individual with a total of only 14
 publications to his credit who spends 7 years sulking in his attic
 would be a strong candidate for early retirement at an aggressive
 British research department.

\subsection*{Fermat--Wiles in three minutes}
{\bf Fermat's Last Theorem:} {\em A perfect cube cannot be written as
the sum of two perfect cubes, a perfect fourth power cannot be written
as the sum of two perfect fourth powers, and likewise, a perfect nth
power cannot be written as the sum of two perfect nth powers. In other
words, for any $n>2$, the equation
 \begin{equation}
 \renewcommand\theequation{$*$}
 a^n + b^n = c^n
 \end{equation}
does not have any integer solutions with $a,b,c\ne0$.} \medskip

Over the 350 years since Fermat's celebrated margin, any number of
mathematicians have tried their hands at this, from 10 year olds in
public libraries through to the most distinguished professors. A
popular approach is to argue by contradiction: if $a,b,c$ are nonzero
integers satisfying Fermat's equation $(*)$, you try to argue that
$a,b,c$ are very special, in fact eventually so special that they can't
exist. Any prime dividing the right-hand side of Fermat's equation $(*)$
divides it $n$ times, and you could try to argue that it can't also
divide the left-hand side $n$ times. About 150 years ago, a number of
people noted that the left-hand side splits as a product of $n$ factors
in the ring of cyclotomic integers, these factors being more-or-less
coprime, and thought that they could see a way through from this to a
contradiction; in the course of explaining why this approach fails,
Kummer invented algebraic number theory and the class group of an
algebraic number field, and paved the way for class field theory.

A key twist on the argument by contradiction was invented in the early
1980s by the German mathematician Gerhard Frey: if $a,b,c$ are nonzero
integers satisfying Fermat's equation $(*)$, consider the equation
 \begin{equation}
 \renewcommand\theequation{$**$}
 y^2 = x(x+a^n)(x+c^n),
 \end{equation}
where $a,b,c$ are considered fixed. This equation in $x,y$ is called an
{\em elliptic curve:} it is the curve obtained as the graph of the
function square root of $x(x+a^n)(x+c^n)$. (The name ``elliptic''
comes from the fact that equations of this form arise in Euler's
integral formula for the arc length of an ellipse.) Just as before, the
aim is to argue that Frey's curve $(**)$ is very special, in fact
eventually so special that it can't exist. (The special thing is that
the discriminant of the cubic polynomial on r-h.s. of $(**)$ is
$a^nb^nc^n$, which has many repeated prime factors.) Frey's idea was
immediately taken up by a number of mathematicians, who hoped to
exploit the encyclopaedia of results on elliptic curves accumulated
since the time of Fermat and Euler.

The deepest fact about elliptic curves, and the essential achievement of
Wiles' work from the mid 1980s, is the Taniyama--Shimura conjecture:
every elliptic curve over the integers is ``modular'', that is,
parametrised by modular forms. A modular form is a function having very
strong symmetry with respect to an arithmetic group -- it is thus an
object of complex analysis, hyperbolic geometry, representation theory,
arithmetic and algebraic geo\-metry. It was known from the late 1980s
that the Taniyama--Shimura conjecture would imply that the Frey curves
do not exist, hence prove Fermat's last theorem. The received wisdom
from the 1970s was that Taniyama--Shimura was most likely to be true,
but unlikely ever to be proved. So it might well have remained without
Wiles' 7 year odyssey in his attic. For the details of the proof, I
borrow a phrase from Fermat: {\em Hoc elogium exiguitas non caperet}.

In one sense, all this talk of Fermat and what's happened over the last
350 years is extremely misleading, because the real impact of Wiles' work
lies in the future. Possibly the single biggest issue in the mathematics
of the next century is a vast generalisation of class field theory and of
the Taniyama--Shimura conjecture called the {\em Langlands program:}
 \begin{quote}
  modular forms parametrise the representations of\\
  the Galois group of the rational number field. 
 \end{quote}
While it certainly sees off Fermat's last theorem, Wiles's work, and
its current development at his hands and those of his students and
successors, is a searchlight illuminating this maze.

\subsection*{Presentation}
Mr Chancellor, in the name of the council, I present to you for
admission to the degree of Doctor of Science, {\em honoris causa},

\bigskip 
\qquad {\bf Andrew Wiles}

\bigskip \noindent Miles Reid, Jul 1998

 \end{document}